\documentclass{amsart} 
 
\usepackage{amsmath,amssymb}

\begin{document}

 %%%%%%%%%%%%%%%%%%%%%%%%
 \newtheorem{theorem}{Theorem}[section]
\newtheorem{lemma}[theorem]{Lemma}
\newtheorem{proposition}[theorem]{Proposition}
\newtheorem{corollary}[theorem]{Corollary}
\newtheorem{remark}[theorem]{Remark}
\newtheorem{example}[theorem]{Example}
\newtheorem{question}[theorem]{Question}
\newtheorem{questions}[theorem]{Questions}
\newtheorem{problem}[theorem]{Problem}
\newtheorem{claim}[theorem]{Claim}
\newtheorem{definition}[theorem]{Definition}
%%%%%%%%%%%%%%%%%%%%%%%

%%%%%%% Preamble

\title{Heredity of $\tau$-pseudocompactness}

\author{Jerry E. Vaughan}

\date{30 April 2005}

\thanks{Primary 54B05, 54D80, 54G20, 54F05.
}

%\address{Department of Mathematical Sciences\\
%UNC-Greensboro\\
%Greensboro, NC, 27402\\
%}

%\email{vaughanj@uncg.edu}

\begin{abstract} S. Garc\'{\i}a-Ferreira and H. Ohta gave a construction that was intended to produce a
$\tau$-pseudocompact space, which has a regular-closed zero set $A$ and a regular-closed $C$-embedded set $B$ such 
that neither $A$ nor $B$ is 
$\tau$-pseudocompact.  We show that although their  sets $A$, $B$  are not regular-closed, there are at least two ways to
make their construction work to give the desired example.
\end{abstract}

%%%End Preamble

\maketitle

\section{Introduction} All spaces considered in this paper are
Tychonoff, i.e., $T_{3\frac{1}{2}}$-spaces. Let $\tau\geq\omega$ denote an infinite cardinal number, 
and $\mathbb R^\tau$ the product of $\tau$ copies of the real line with the product topology.   
J. F. Kennison defined a space $X$ to be {\em $\tau$-pseudocompact} 
provided for every continuous $f:X\rightarrow
\mathbb R^\tau$, $f(X)$ is a closed subset of $\mathbb R^\tau$ \cite{kennison}. He proved that a space $X$ is
$\tau$-pseudocompact  if and only if whenever  $\mathcal F$ is a family
of zero sets of
$X$ with the finite intersection property (FIP) and $|\mathcal F|\leq\tau$, then $\cap\mathcal
F\not=\emptyset$ \cite[Theorem 2.2]{kennison}. It is known and easy to prove that 
$\omega$-pseudocompactness is equivalent to the well-known notion of pseudocompactness (e.g., see
\cite[Theorem 2.1]{kennison}).

Recall that a subset $H$ of a topological space is called {\em regular-closed} if $H$ is
the closure of an open set.  $H$ is called a {\em zero set} provided there exists a continuous
$f:X\rightarrow [0,1]$ such that $H= f^{-1}(0)$,  and $H$ is called {\em $C$-embedded in $X$} if for every continuous 
$f:H\rightarrow\mathbb R$, there is a continuous  $g:X\rightarrow\mathbb R$ such that $g$
extends $f$. A set $Y\subset X$ is said to be {\em countably compact in $X$} if every infinite subset of $Y$ has a limit
point in
$X$ \cite{gfo}

There are several  known examples  that show $\tau$-pseudocompactness is not hereditary to various kinds
of closed sets.  Kennison showed that
$\tau$-pseudocompactness is not hereditary to closed $C$-embedded sets
\cite[p.440]{kennison}.  T. Retta showed  (for
$\tau\geq
\mathfrak c$) that
$\tau$-pseudocompactness is not hereditary to regular-closed subsets \cite{retta}, and a different
construction to show the same thing was given by S. Garc\'{\i}a-Ferreira, M. Sanchis, and S. Watson
\cite[Corollary 1.4]{gsw}, assuming $cf(\tau)>2^{ \mathfrak c}$.  These examples 
demonstrate a difference between the countable and uncountable cases: pseudocompactness (i.e.,
$\omega$-pseudocompactness) is hereditary to $C$-embedded subsets and to regular-closed sets (e.g., see 
\cite[9.13]{GJ}), but for $\tau\geq \mathfrak c$, $\tau$-pseudocompactness is not necessarily hereditary to either kind of
closed set. Concerning cardinals not covered by the previous examples,
 H. Ohta (see \cite{gsw}) constructed an example  to show that $\omega_1$-pseudocompactness is not hereditary
to regular-closed sets.  Garc\'{\i}a-Ferreira and  Ohta
\cite[Example 2.4]{gfo} generalized this construction to all uncountable cardinals. They  stated the following

\begin{example}[Garc\'{\i}a-Ferreira and  Ohta]\label{main} For all $\tau\geq\omega_1$, there exists
a $\tau$-pseudocompact space 
$X$ with two regular-closed sets
$A,B$ such that
$A$ is a zero set, $B$ is $C$-embedded, and neither is $\tau$-pseudocompact.
\end{example}

\medskip
There is, however, a small gap in the constructions of Ohta in \cite{gsw} and of Garc\'{\i}a-Ferreira and  Ohta in
\cite{gfo}.  The purpose of this paper is to show in
\S 2 that the sets $A$ and $B$ that they claim in \cite{gfo} to be regular-closed are not, and to show in
\S 3  that a simple modification of their construction suffices to prove Example \ref{main}. The
modification is to replace the cardinals
$\tau^+$ and
$\omega_1$ in the Garc\'{\i}a-Ferreira and  Ohta construction with their long line counterparts. Possibly the previous
sentence is sufficient for our main goal of establishing Example \ref{main}, but we elaborate a bit more on this in
$\S 2$. In \S 3 we present another way to modify their construction and give a different, possibly simpler, proof of Example
\ref{main}. 
\medskip

Garc\'{\i}a-Ferreira and  Ohta  also proved that $\tau$-pseudocompactness is hereditary to any subset that
is both a zero set and a $C$-embedded set (regular-closed or not) \cite[Theorem. 1.4]{gfo}. Thus Example \ref{main} seems 
to be about as strong as possible, and is therefore an important example in the theory of $\tau$-pseudocompactness.

\section{The Construction of Garc\'{\i}a-Ferreira and Ohta}\label{gf-o}

First we recall the Alexandroff duplicate $A(X)$ of a space $X$.  The underlying set of $A(X)$ is $X\times 2$,
where
$2 = \{0,1\}$. In the topology of $A(X)$, each point of $X\times\{1\}$ is isolated, and each point $(x,0)\in X\times\{0\}$
has basic open neighborhoods of the form $U\times 2\setminus\{(x,1)\}$, where $U$ is an open neighborhood of $x$ in $X$ 
(see
\cite{e}).  Let $Y\subset X$. The space $A(X,Y)$ is defined to be the set $(X\times\{0\})\cup (Y\times\{1\})$ with the
subspace topology from
$A(X)$
\cite[\S 2]{gfo}.

Now we recall the construction of Garc\'{\i}a-Ferreira and Ohta \cite[Example 2.4]{gfo}.  
Let $\tau\geq\omega_1$
be an infinite  cardinal, and $\tau^+$ the first cardinal larger than $\tau$.  As is well known,
the spaces
$\tau^+$ and $\omega_1$  with the order topology satisfy the following properties:
\smallskip

(1) $\tau^+$  is initially
$\tau$-compact (i.e., every open cover of cardinality at most $\tau$ has a finite subcover \cite{AU}) and 
$\omega_1$ is initially
$\omega$-compact (i.e., countably compact).  

\smallskip
(2)  every real-valued continuous function defined on $\tau^+$ or $\omega_1$ is eventually constant.
\smallskip

\medskip
 Let
$S_1 = (\tau^++1)\times(\omega+1)$, and $S_2 = 
(\omega_1+1)\times(\omega+1)$.  Next consider the quotient of the
disjoint union $S_1\oplus S_2$ obtained by identifying $(\tau^+,n)$ and $(\omega_1,n)$ for every
$n\in\omega$.  Let $\varphi$ denote the quotient map from $S_1\oplus S_2$ onto the quotient space.  Then let $X$ denote
the quotient space minus the point $\varphi((\tau^+,\omega))=\varphi((\omega_1,\omega))$.  Let $Y_1 = \tau^+\times\{\omega\}\subset S_1$,
 $Y_2= \omega_1\times\{\omega\}\subset S_2$, and $Y =\varphi(Y_1\cup Y_2)$, and $Z =\varphi(Y_2)$.

The space for Example \ref{main} given by Garc\'{\i}a-Ferreira and  Ohta is $A(X,Y)$ where $X,Y$ were defined
in the previous paragraph, and the two subsets are $A= Y\times 2$ and $B =Z\times 2$.

A gap in the proof by Garc\'{\i}a-Ferreira and  Ohta occurs because neither of $A= Y\times 2$ or $B =Z\times 2$ is
a regular-closed set.  To see this let $int_X(H)$ denote the interior of $H$ in $X$, and note the following fact: For  
any space
$X$, if 
$H$ is closed in $X$ and there is a point $p
\in H$ such that $p\not\in int_X(H)$ and $p$ is relatively isolated in $H$, then
$H$ is not regular-closed.  For the sets
$A= Y\times 2$ and $B =Z\times 2$, take any isolated ordinal $\alpha <\omega_1$ and put $p = (\varphi(\alpha),0)$.
Then  $p$ is relatively isolated in $A$ and $B$, hence $A$, $B$ are not regular-closed.

\medskip
The following Lemma indicates a way to repair this gap.

\begin{lemma}\label{regularclosed} If $Y\subset X$ is dense-in-itself, then $Y\times 2$ is regular-closed in $A(X,Y)$.
\end{lemma}
Proof.  We claim that $Y\times 2 = cl_{A(X,Y)}(Y\times\{1\})$.  Since $Y\times 2$ is closed in $A(X,Y)$, we need only  show
that
 $Y\times\{0\}\subset cl_{A(X,Y)}(Y\times \{1\})$. For any $y\in Y$, and any neighborhood $U$ of $y$ in $X$, $(U\times
2)\setminus\{(y,1)\}$ is a basic  neighborhood of the point $(y,0)$ in $A(X)$. Since $Y$ is dense-in-itself there is
$z\in U\cap Y$ such that $z \not= y$. Then $(z,1) \in (U\times
2)\setminus\{(y,1)\}$ which shows that $(y,0)\in cl_{A(X,Y)}(Y\times \{1\})$.

\section{The First Modification}
To repair Example \ref{main} we start over the construction of Garc\'{\i}a-Ferreira and  Ohta, but this time we use the long
line counterparts of the cardinals $\tau^+$, and $\omega_1$.  The following lemmas indicate that the counterparts have the
key properties needed in the construction, and since each of these counterparts is dense-in-itself (in fact, connected), 
Lemma \ref{regularclosed} fixes the gap and Example \ref{main}  follows.  

\medskip
Notation: Fix an uncountable cardinal $\tau$.  Let   $T =
\tau^+\times_{lex}[0,1)$ and $W = \omega_1\times_{lex} [0,1)$ where the products are given the lexicographic order and
the order topology.

\begin{lemma} $W$ is countably compact, and $T$ is initially $\tau$-compact.
\end{lemma}

\begin{lemma} (cf.  \cite[16H]{GJ}) Every real-valued  continuous function defined on $W$ or $T$ is eventually constant.
\end{lemma}

To get counterparts to $\tau^++1$ and $\omega_1+1$, let  $W+1=W\cup \{w\}$ and $T+1=T\cup\{t\}$, where $w,t$ are points
not in $T\cup W$.  Extend the order of $W$ and $T$ so that $w$ acts as the last element of $W$ and $t$ acts as the last
element of $T$.
 Let
$S_1 = (T+1)\times(\omega+1)$, and $S_2 = 
(W+1)\times(\omega+1)$.  Next consider the quotient of the
disjoint union $S_1\oplus S_2$ obtained by identifying $(t,n)$ and $(w,n)$ for every
$n\in\omega$.  Let $\varphi$ denote the quotient map from $S_1\oplus S_2$ onto the quotient space.  Then let $X$ denote
the quotient space minus the point $\varphi((t,\omega))=\varphi((w,\omega))$.  Let $Y_1 =T\times\{\omega\}\subset S_1$,
 $Y_2= W\times\{\omega\}\subset S_2$, and $Y =\varphi(Y_1\cup Y_2)$, and $Z =\varphi(Y_2)$.  

The space for Example \ref{main} is $A(X,Y)$ where $X,Y$ were defined
in the previous paragraph, and the two subsets are $A= Y\times 2$ and $B =Z\times 2$.  Since
$T$ and $W$ are connected (e.g., see \cite[16H]{GJ}),  each of the sets $Y$ and $Z$ is dense-in-itself, so by Lemma
\ref{regularclosed} they are regular-closed in $A(X,Y)$. The other properties required in Example \ref{main} follow as in
\cite{gfo}.

\section{Another Modification}\label{second}

In this section we present another  modification of the construction of
Garc\'{\i}a-Ferreira and  Ohta, suggested to us by Alan Dow,  which gives a second proof of Example \ref{main}. This \
modification does not use lexicographic products. First we formalize a variation
of the Alexandroff duplicate  construction, which is probably not new.

Let $C=\{\frac{1}{n}: n\geq 1\}\cup\{0\}$ denote the usual convergent sequence.  Let $X$ be a space and put
$M(X) = X\times C$. Define a topology on $M(X)$ as follows.  All points of the form $(x,\frac{1}{n}) $ for $n\geq 1$
are isolated, and basic neighborhoods for a point $(x,0)$ are defined to be sets of the form
$(U\times C)\setminus F$ where $U$ is an open neighborhood of $x$ in $X$ and $F$ is a finite set.  It is routine to check
that this topology on $M(X)$ is $T_{3\frac{1}{2}}$.

For $Y\subset X$, we define $M(X,Y) = (X\times\{0\})\cup (Y\times\{\frac{1}{n}:n\geq 1\})$ with the subspace topology from
$M(X)$. Note that $M(X)=M(X,X)$.  Let $\pi$ denote the projection
map $\pi:M(X)\rightarrow X$ defined by  $\pi(x,e) = x$ for all $e\in C$. By abuse of notation we also let $\pi$ denote
the restriction of this projection map to $M(X,Y)$.

\begin{lemma}If $Y$ is a zero set of $X$, then $Y\times C$ is a zero set of $M(X,Y)$.
\end{lemma}

Proof. This follows because the projection map $\pi$ is continuous.

\begin{lemma}
If $Z\subset Y$ is $C$-embedded in $X$ then $Z\times C$ is $C$-embedded in $M(X,Y)$. 
\end{lemma}
Proof. Given a continuous function $f:(Z\times C)\rightarrow \mathbb R$, we may continuously extend 
$f\upharpoonright(Z\times \{0\})$ to all of $X\times\{0\}$ because $Z$ is $C$-embedded in $X$; so we may assume $f$ is defined
on
$X\times\{0\}\cup Z\times C$.    Then define
$g:M(X,Y)\rightarrow
\mathbb R$ by

 \[g(p) = \left\{ \begin{array}{ll}
f(p) & \mbox{if $p\in X\times\{0\}\cup Y\times C$ }\\
f((y,0)) & \mbox{if $p=(y,e)$ and $y\in Y\setminus Z$}
\end{array}\right.\]
The function $g$ is continuous by a standard gluing lemma, and extends $f$ to $M(X,Y)$.
\medskip

The next result is an analog of \cite[Lemma 2.2]{gfo}.

\begin{lemma}\label{analog} M(X,Y) is $\tau$-pseudocompact if and only if $X$ is $\tau$-pseudocompact and $Y$ is
countably compact in $X$. 
\end{lemma}

Proof.  Assume that  $X$ is $\tau$-pseudocompact and $Y$ is countably compact in
$X$.  Let $\{f_\beta^{-1}(0):\beta <\tau\}$ be a family of zero sets of $M(X,Y)$ with the FIP.
If this family traces on $X\times\{0\}$, then the intersection is non-empty because $X$ is $\tau$-pseudocompact.
Thus we
assume there is
$\alpha<\tau$ such that $f_\alpha^{-1}(0)\subset Y\times\{\frac{1}{n}:n\geq 1\}$.  Note that if $H\subset 
Y\times\{\frac{1}{n}:n\geq 1\}$ and $H$ is  closed in $M(X,Y)$ then $H$  is finite. This is because either 
$\{y\in Y:(\exists n\geq 1)((y,\frac{1}{n})\in H)\}$ is infinite, hence
  has a limit point
$x\in X$ which implies
$(x,0)\in cl_{M(X,Y)}(H)=H$, or there is a $y\in Y$ such that $(y,\frac{1}{n})\in H$ for infinitely many $n$, hence
$(y,0)\in cl_{M(X,Y)}(H)=H$, which is again a contradiction.  Thus $f_\alpha^{-1}(0)$ is finite, hence compact; so
one of the points in $f_\alpha^{-1}(0)$ is in $\cap\{f^{-1}_\beta(0):\beta<\tau\}$.  Thus $M(X,Y)$ is $\tau$-pseudocompact.

Conversely, suppose $M(X,Y)$ is $\tau$-pseudocompact.  Let  $\mathcal F=\{f_\beta^{-1}(0): \beta<\tau\}$
be a family of zero sets of $X$ with the FIP. Since the projection
map
$\pi$ is continuous, the maps $g_\alpha = f_\alpha\circ \pi$ are continuous on $M(X)$ for all
$\alpha<\tau$. Thus
$\{g_\beta^{-1}(0):\beta<\tau\}$ is a family of zero sets on $M(X)$.
 Since
$$f_\alpha^{-1}(0)\times\{0\}\subset g_\alpha^{-1}(0)\cap (X\times\{0\})\subset M(X,Y),$$
$\mathcal G =\{g^{-1}(0)\cap M(X,Y):\alpha<\tau\}$ has the FIP.  By assumption, there exists $p\in\cap\mathcal G$; so
$\pi(p)\in\cap\mathcal F$. Thus $X$ is $\tau$-pseudocompact. To see that $Y$ is countably compact in $X$, suppose
otherwise, i.e., suppose  there is an infinite subset of
$Y$ that  has no limit points   in $X$.  Then there is an infinite subset of $Y\times\{1\}$ that forms a closed discrete set
of isolated points, but this is impossible because $M(X,Y)$ is pseudocompact.

\medskip

To complete the construction,
 let $X$, $Y$, $Z$ be the space and subsets defined in \S \ref{gf-o} and put $A = Y\times C$, and $B= Z\times C$.  Clearly
$B$ is clopen in $A$.   Since $Y$ is  a zero set in $X$, 
$A = Y\times C$ is a zero set in $M(X,Y)$.  Further  $A=cl_{M(X,Y)}(Y\times\{\frac{1}{n}:n\geq 1\})$; so $A$ is
a regular-closed set (although $Y$ is not dense-in-itself). Similarly, $B$ is regular-closed.  

\medskip
To complete our second proof of Example \ref{main}, we use the next lemma 
which follows the method of Garc\'{\i}a-Ferreira and  Ohta.

\begin{lemma} $A$ and $B$ are not $\omega_1$-pseudocompact.
\end{lemma}
Proof. Since $B$ is a clopen subset of $A$, it suffices to show that $B$ is not $\omega_1$-pseudocompact.
Now $B=Z\times C =
M(Z)=M(Z,Z)$. Since
$Z$ is a copy of
$\omega_1$,
$Z$  contains a decreasing   family
 of $\omega_1$ many clopen sets with empty intersection; so $Z$ is not $\omega_1$-pseudocompact. By Lemma~\ref{analog},
 $M(Z) = B$ is not $\omega_1$-pseudocompact, hence $A$ is not.

\medskip
We thank Alan Dow for  suggesting the space discussed in \S \ref{second} and for helpful remarks concerning it.
%%%%%%Bibliography

\end{document}